\documentclass[secthm,dvips]{arxbj}
\usepackage{graphicx}

\usepackage{accents}

\aid{0}
\volume{16}
\issue{3}
\pubyear{2010}
\firstpage{705}
\lastpage{729}
\doi{10.3150/09-BEJ228}

\makeatletter
\newproclaim{definition}{Definition}[section]
\newtheorem{proposition}[definition]{Proposition}
\newtheorem{propositionn}[thm]{Proposition}
\newproclaim{condition}{Condition}
\newremark{example}[thm]{Example}
\newtheorem{lemma}[thm]{Lemma}
\makeatother

\begin{document}
\begin{frontmatter}

\title{Functional linear regression via canonical analysis}
\runtitle{Regression via canonical analysis}

\begin{aug}
\author[1]{\fnms{Guozhong} \snm{He}\thanksref{1}\ead[label=e1]{ghe@dhs.ca.gov}},
\author[2]{\fnms{Hans-Georg} \snm{M\"{u}ller}\corref{}\thanksref{2,e2}\ead[label=e2,mark]{mueller@wald.ucdavis.edu}},
\author[2]{\fnms{Jane-Ling} \snm{Wang}\thanksref{2,e3}\ead[label=e3,mark]{wang@wald.ucdavis.edu}}\\
\and
\author[2]{\fnms{Wenjing} \snm{Yang}\thanksref{2,e4}\ead[label=e4,mark]{wyang@wald.ucdavis.edu}}
\runauthor{He, M\"{u}ller, Wang and Yang}
\address[1]{University of California, San Francisco, Institute for Health and Aging and
California Diabetes Program, California Department of Public Health P.O. Box
997377, MS 7211, Sacramento, CA 95899-7377, USA. \printead{e1}}
\address[2]{Department of Statistics, University of California, One Shields
Avenue, Davis, CA 95616,
USA.\\
\printead{e2}; \printead*{e3}; \printead*{e4}}
\end{aug}

\received{\smonth{11} \syear{2007}}
\revised{\smonth{6} \syear{2009}}

%
\begin{abstract}
We study regression models for the situation where both dependent
and independent variables are square-integrable stochastic
processes. Questions concerning the definition and existence of the
corresponding functional linear regression models and some basic
properties are explored for this situation. We derive
a~representation of the regression parameter function in terms of the
canonical components of the processes involved. This representation
establishes a connection between functional regression and
functional canonical analysis and suggests alternative approaches
for the implementation of functional linear regression analysis. A
specific procedure for the estimation of the regression parameter
function using canonical expansions is proposed and compared with an
established functional principal component regression approach. As
an example of an application, we present an analysis of mortality data for
cohorts of medflies, obtained in experimental studies of aging and
longevity.
\end{abstract}

%
\begin{keyword}
\kwd{canonical components}
\kwd{covariance operator}
\kwd{functional data analysis}
\kwd{functional linear model}
\kwd{longitudinal data}
\kwd{parameter function}
\kwd{stochastic process}
\end{keyword}

\end{frontmatter}

\section{Introduction}

With the advancement of modern technology, data sets which
contain repeated measurements obtained on a dense grid are becoming
ubiquitous. Such data can be viewed as a sample of curves
or functions and are referred to as functional data. We consider
here the extension of the linear regression model to the case of
functional data. In this extension, both predictors and responses
are random functions rather than random vectors. It is well known
(Ramsay and Dalzell (\citeyear{Ramsay91}); Ramsay and Silverman (\citeyear{Ramsay05})) that the
traditional linear regression model for multivariate data, defined
as
%
%
\begin{equation}\label{basic}
\mathbf{Y}=\bolds{\alpha}_{0}+\mathbf{X}\bolds{\beta
}_{0}+\bolds{
\varepsilon},
\end{equation}
may be extended to the functional setting by postulating the model,
for $s\in T_{1},t\in T_{2}$,
%
%
\begin{equation}\label{linear}
Y(t)=\alpha_{0}(t)+\int_{T_{1}}X(s) \beta_{0}(s,t)\,\mathrm{d}s+ \varepsilon
(t).
\end{equation}

Writing all vectors as row vectors in the classical model
(\ref{basic}), $\mathbf{Y}$ and $\bolds{\varepsilon}$ are random vectors
in~$\mathbf{R}^{p_{2}}$, $\mathbf{X}$ is a random vector in
$\mathbf{R}^{p_{1}}$, and $\bolds{\alpha}_{0}$ and
$\bolds{\beta}_{0}$ are, respectively, $1\times p_{2}$ and $p_{1}\times p_{2}$ matrices
containing the regression parameters. The vector $\bolds{\varepsilon}$
has the usual interpretation of an error vector, with
$E[\bolds{\varepsilon}]=0$ and $\operatorname{cov}[\bolds
{\varepsilon}]=\sigma^{2}I$,
$I$ denoting the identity matrix. In the functional model
(\ref{linear}), random vectors $\mathbf{X}, \mathbf{Y}$ and $\bolds
{\varepsilon}$ in (\ref{basic}) are replaced by random functions defined on the
intervals~$T_{1}$ and $T_{2}$. The extension of the classical linear
model (\ref{basic}) to the functional linear model (\ref{linear}) is
obtained by replacing the matrix operation on the right-hand side of
(\ref{basic}) with an integral operator in (\ref{linear}). In the
original approach of Ramsay and Dalzell (\citeyear{Ramsay91}), a penalized
least-squares approach using L-splines was adopted and applied to a study
in temperature-precipitation patterns, based on data from Canadian
weather stations.

The functional regression model (\ref{linear}) for the case of
scalar responses has attracted much recent interest (Cardot and Sarda
(\citeyear{Cardot05}); M\"{u}ller and Stadtm\"{u}ller (\citeyear{Muller05}); Hall and Horowitz (\citeyear{Hall07})),
while the case of functional responses has been much less thoroughly
investigated (Ramsay and Dalzell (\citeyear{Ramsay91}); Yao, M\"{u}ller and Wang
(\citeyear{Yao05b})). Discussions on various approaches and estimation procedures
can be found in the insightful monograph of Ramsay and Silverman
(\citeyear{Ramsay05}). In this paper, we propose an alternative approach to predict
$Y(\cdot)$ from $ X(\cdot)$, by adopting a novel canonical
representation of the regression parameter function $\beta
_{0}(s,t)$. Several distinctive features of functional linear models
emerge in the development of this canonical expansion approach.

It is well known that in the classical multivariate linear model,
the regression slope parameter matrix is uniquely determined by
$\bolds{\beta}_{0}=\operatorname{cov}(\mathbf{X})^{-1}\operatorname
{cov}(\mathbf{X,Y})$, as long as the
covariance matrix $\operatorname{cov}(\mathbf{X})$ is invertible. In
contrast, the
corresponding parameter function $\beta_{0}(\cdot,\cdot)$,
appearing in (\ref{linear}), is typically not identifiable. This
identifiability issue is discussed in Section \ref{sec2}. It relates to the
compactness of the covariance operator of the process $X$ which
makes it non-invertible. In Section \ref{sec2}, we demonstrate how
restriction to a subspace allows this problem to be circumvented. Under
suitable restrictions, the components of model (\ref{linear}) are
then well defined.

Utilizing the canonical decomposition in Theorem \ref{th3.3} below leads to
an alternative approach to estimating the parameter function $\beta
_{0}(\cdot,\cdot)$. The canonical decomposition links $Y$ and $X$
through their functional canonical correlation structure. The
corresponding canonical components form a bridge between canonical
analysis and linear regression modeling. Canonical components
provide a decomposition of the structure of the dependency between
$Y$ and $X$ and lead to a natural expansion of the regression
parameter function $\beta_{0}(\cdot,\cdot)$, thus aiding in its
interpretation. The canonical regression decomposition also suggests
a new family of estimation procedures for functional regression
analysis. We refer to this methodology as \textit{functional canonical
regression analysis}. Classical canonical correlation analysis (CCA)
was introduced by Hotelling (\citeyear{Hotelling36}) and was connected to function
spaces by Hannan (\citeyear{Hannan61}). Substantial extensions and connections to
reproducing kernel Hilbert spaces were recently developed in Eubank
and Hsing (\citeyear{Eubank08}); for other recent developments see Cupidon \textit{et
al.}~(\citeyear{Cupidon07}).

Canonical
correlation is known not to work particularly well for very
high-dimensional multivariate data, as it involves an inverse problem.
Leurgans, Moyeed and Silverman (\citeyear{Leurgans93}) tackled the difficult problem
of extending CCA to the case of infinite-dimensional functional data
and discussed the precarious regularization issues which are faced; He,
M\"{u}ller and Wang (\citeyear{He03,He04}) further explored various aspects
and proposed practically feasible regularization procedures for
functional CCA. While CCA for functional data is worthwhile, but
difficult to implement and interpret, the canonical approach to
functional regression is here found to compare favorably with the
well established principal-component-based regression approach in an
example of an application (Section \ref{sec5}). This demonstrates a potentially
important new role for canonical decompositions in functional
regression analysis. The functional linear model (\ref{linear})
includes the varying coefficient linear model studied in Hoover \textit{et
al.}~(\citeyear{Hoover98}) and Fan and Zhang (\citeyear{Fan98}) as a special case, where $\beta
(s,t)=\beta(t)\delta_{t}(s)$; here, $\delta_{t}(\cdot)$ is a
delta function centered at $t$ and~$\beta(t)$ is the varying
coefficient function. Other forms of functional regression models
with vector-valued predictors and functional responses were
considered by Faraway (\citeyear{Faraway97}), Shi, Weiss and Taylor (\citeyear{Shi96}), Rice and
Wu (\citeyear{Rice00}), Chiou, M\"{u}ller and Wang (\citeyear{Chiou03}) and Ritz and Streibig
(\citeyear{Ritz09}).

The paper is organized as follows. Functional canonical analysis and
functional linear models for $L_{2}$-processes are introduced in
Section \ref{sec2}. Sufficient conditions for the existence of functional
normal equations are given in Proposition \ref{pr2.2}. The canonical
regression decomposition and its properties are the theme of Section
\ref{sec3}. In Section \ref{sec4}, we propose a novel estimation technique to obtain
regression parameter function estimates based on functional
canonical components. The regression parameter function is the basic
model component of interest in functional linear models, in analogy
to the parameter vector in classical linear models. The proposed
estimation method, based on a canonical regression decomposition, is
contrasted with an established functional regression method based on
a principal component decomposition. These methods utilize a
dimension reduction step to regularize the solution of the inverse
problems posed by both functional regression and functional
canonical analysis. As a selection criterion for tuning parameters,
such as bandwidths or numbers of canonical components, we use
minimization of prediction error via leave-one-curve-out
cross-validation (Rice and Silverman (\citeyear{Rice91})). The proposed estimation
procedures are applied to mortality data obtained for cohorts of
medflies (Section \ref{sec5}). Our goal in this application is to predict a
random trajectory of mortality for a female cohort of flies from the
trajectory of mortality for a male cohort which was raised in the
same cage. We find that the proposed functional canonical regression
method gains an advantage over functional principal component
regression in terms of prediction error.

Additional results on canonical regression decompositions and
properties of functional regression operators are compiled in
Section \ref{sec6}. All proofs are collected in Section \ref{sec7}.

\section{Functional linear regression and the functional
normal equation}\label{sec2}

In this section, we explore the formal setting as well as
identifiability issues for functional linear regression models. Both
response and predictor functions are considered to come from a
sample of pairs of random curves. A basic assumption is that all
random curves or functions are square-integrable stochastic
processes. Consider a measure $\mu$ on a real index set $T$ and
let $L_{2}(T)$ be the class of real-valued functions such that
$\|f\|^{2}=\int_{T}|f|^{2}\,\mathrm{d}\mu<\infty$. This is a Hilbert space
with the inner product $\langle f,g \rangle=\int_{T}fg\,\mathrm{d}\mu$ and we
write $f=g$ if $\int_{T}(f-g)^{2}\,\mathrm{d}\mu=0$. The index set $T$ can be
a set of time points, such as $T=\{1,2,\ldots,k\}$, a compact
interval $ T=[a,b]$ or even a rectangle formed by two intervals
$S_{1}$ and $S_{2}$, $T=S_{1}\times S_{2}$. We focus on index sets
$T$ that are either compact real intervals or compact rectangles in
$\mathbf{R}^{2}$ and
consider $\mu$ to be the Lebesgue measure on~$\mathbf{R}^{1}$ or
$\mathbf{R}%
^{2}$. Extensions to other index sets $T$ and other measures are
self-evident. An $L_{2}$-process is a stochastic process
$X=\{X(t),t\in T\}$, $X\in L_{2}(T)$, with $E[\|X\|^{2}]< \infty,
E[X(t)^2]<\infty$ for all $t \in T$.

Let $X\in L_{2}(T_{1})$ and $Y\in L_{2}(T_{2})$.

\begin{definition} Processes $(X,Y)$ are
subject to a
functional linear model if
%
%
\begin{equation}\label{linear1}
Y(t)=\alpha_{0}(t)+\int_{T_{1}}X(s)\beta_{0}(s,t)\,\mathrm{d}s+\varepsilon
(t),\qquad t \in T_2,
\end{equation}
where $\beta_{0}\in L_{2}(T_{1}\times T_{2})$ is the
parameter function, $\varepsilon\in L_{2}(T_{2})$ is a random error
process with $E[\varepsilon(t)]=0$ for $t \in T_1$, and
$\varepsilon$ and $X$ are uncorrelated, in the sense that
$E[X(t)\varepsilon(s)]=0$ for all $s,t \in T_1$.
\end{definition}

Without loss of generality, we assume from now on that all processes
considered have zero mean functions, $EX(t)=0$ and $EY(s)=0$ for all
$t$, $s$. We define the regression integral operator
$\mathcal{L}_X: L_{2}(T_{1}\times T_{2})\rightarrow
L_{2}(T_{2})$ by
\[
(\mathcal{L}_X\beta)(t)=\int_{T_{1}}X(s)\beta(s,t)\,\mathrm{d}s
\qquad\mbox{for any } \beta\in L_{2}(T_{1}\times T_{2}).
\]
Equation (\ref{linear1}) can then be rewritten as
%
%
\begin{equation}\label{linear2}
Y=\mathcal{L}_X\beta_{0}+\varepsilon.
\end{equation}
Denote the auto- and cross-covariance functions of $X$ and $Y$ by
\begin{eqnarray*}
r_{XX}(s,t)&=&\operatorname{cov}[X(s),X(t)],\qquad s, t\in
T_{1},\\
r_{YY}(s,t)&=&\operatorname{cov}[Y(s),Y(t)],\qquad s,t \in T_{2} ,
\quad\mbox{and}\\
r_{XY}(s,t)&=&\operatorname{cov}[X(s),Y(t)],\qquad s\in T_{1}, t\in T_{2}.
\end{eqnarray*}
The
autocovariance operator of $X$ is the integral operator
$R_{XX}\dvtx L_{2}(T_{1})\rightarrow L_{2}(T_{1})$, defined by
\[
(R_{XX}u)(s)=\int_{T_{1}}r_{XX}(s,t)u(t)\,\mathrm{d}t,\qquad u\in
L_{2}(T_{1}).
\]
Replacing $r_{XX}$ by $r_{YY}$, $r_{XY}$, we analogously define
operators $R_{YY}\dvtx L_{2}(T_{2})\rightarrow L_{2}(T_{2})$ and
$R_{XY}\dvtx L_{2}(T_{2})\rightarrow L_{2}(T_{1})$, similarly
$R_{YX}$. Then $R_{XX}$ and $R_{YY}$ are compact, self-adjoint and
non-negative definite operators, and $R_{XY}$ and $R_{YX}$ are
compact operators (Conway (\citeyear{Conway85})). We refer to He \textit{et al.}~(\citeyear{He03})
for a
discussion of various properties of these operators.

Another linear operator of interest is the integral operator $\Gamma
_{XX}\dvtx L_{2}(T_{1}\times T_{2})\rightarrow L_{2}(T_{1}\times T_{2})$,
%
%
\begin{equation}\label{gamma}
(\Gamma_{XX}\beta)(s,t)=\int_{T_{1}}r_{XX}(s,w)\beta(w,t) \,\mathrm{d}w.
\end{equation}
The operator
equation
%
%
\begin{equation}\label{norm}
r_{XY}=\Gamma_{XX}\beta,\qquad \beta\in L_{2}(T_{1}\times T_{2})
\end{equation}
is a direct extension of the least-squares normal equation
and may be referred to as the functional population normal equation.

\begin{proposition}\label{pr2.2}
The following statements are
equivalent for a function $\beta_{0} \in L_{2}(T_{1}\times
T_{2})$:
\begin{longlist}[(a)]
\item[(a)] $\beta_{0}$ satisfies the linear model (\ref{linear2});
\item[(b)] $\beta_{0}$ is a solution of the functional normal
equation (\ref{norm});
\item[(c)] $\beta_{0}$ minimizes $E\|Y-\mathcal{L}_X\beta
\|^{2}$ among all $\beta\in L_{2}(T_{1}\times T_{2})$.
\end{longlist}
\end{proposition}

The proof is found Section \ref{sec7}. In the infinite-dimensional case, the
operator $\Gamma_{XX}$ is a Hilbert--Schmidt operator in the Hilbert
space $L_{2}$, according to Proposition \ref{pr6.6} below. A problem we face
is that it is known from functional analysis that a bounded inverse
does not exist for such operators. A consequence is that the
parameter function $\beta_{0}$ in (\ref{linear1}), (\ref{linear2})
is not identifiable without additional constraints. In a situation
where the inverse of the covariance matrix does not exist in the
multivariate case, a unique solution of the normal equation always
exists within the column space of $\operatorname{cov}(\mathbf{X})$
and this solution
then minimizes $E\|Y-\mathcal{L}_X\beta\|^{2}$ on that space. Our
idea to get around the non-invertibility issue in the functional
infinite-dimensional case is to extend this approach for the
non-invertible multivariate case to the functional case. Indeed, as
is demonstrated in Theorem \ref{th2.3} below, under the additional Condition
\ref{condc1}, the solution of (\ref{norm}) exists in the subspace defined by
the range of $\Gamma_{XX}$. This unique solution indeed minimizes
$E\|Y-\mathcal{L}_X\beta\|^{2}$.

We will make use of the Karhunen--Lo\`{e}ve decompositions (Ash and
Gardner (\citeyear{Ash75})) for $L_{2}$-processes $X$ and $Y$,
%
%
\begin{equation}\label{kl}
X(s)=\sum_{m=1}^{\infty}\xi_{m}\theta_{m}(s),\qquad s\in T_{1}
\quad \mbox{and}\quad  Y(t)=\sum_{j=1}^{\infty}\zeta
_{j}\varphi_{j}(t),\qquad t\in T_{2},
\end{equation}
with random variables $\xi_{m}$, $\zeta_{j}$, $m,j\geq1$, and
orthonormal families of $L_{2}$%
-functions $\{\theta_{m}\}_{m\geq1}$ and $\{\varphi_{j}\}_{j\geq
1}$. Here, $E\xi_{m}=E\zeta_{j}=0$, $E\xi_{m}\xi_{p}=\lambda
_{Xm}\delta_{mp}$, $E\zeta_{j}\zeta_{p}=\lambda_{Yj}\delta
_{jp}$ and $\{(\lambda_{Xm},\theta_{m})\},\{(\lambda
_{Yj}, \varphi_{j})\}$ are the eigenvalues and
eigenfunctions of the covariance operators $R_{XX}$ and $%
R_{YY}$, respectively, with $\sum_{m}\lambda_{Xm}<\infty$,
$\sum_{j}\lambda_{Yj}<\infty$. Note that $\delta_{mj}$ is the
Kronecker symbol with $\delta_{mj}=1$ for $m=j$, $\delta_{mj}=0$
for $m\neq j$.

We consider a subset of $L_{2}$ on which inverses of the operator $%
\Gamma_{XX}$ can be defined. As a~Hilbert--Schmidt operator, $%
\Gamma_{XX}$ is compact and therefore not invertible on $L_{2}.$
According to Conway (\citeyear{Conway85}), page 50, the range of $\Gamma_{XX},$
\[
G_{XX}=\{\Gamma_{XX}h\dvtx h\in L_{2}(T_{1}\times T_{2})\},
\]
is characterized by
%
%
\begin{equation}\label{range}
G_{XX}=\Biggl\{g\in L_{2}(T_{1}\times T_{2})\dvtx \sum_{m,j=1}^{\infty
}\lambda_{Xm}^{-2}|\langle g,\theta_{m}\varphi_{j}\rangle
|^{2}<\infty
, g\perp\mathrm{ ker}(\Gamma_{XX})\Biggr\},
\end{equation}
where $\operatorname{ker}(\Gamma_{XX})=\{h:\Gamma_{XX}h=0\}.$
Defining
\[
G_{XX}^{-1}=\Biggl\{h\in L_{2}(T_{1}\times
T_{2})\dvtx\ h=\sum_{m,j=1}^{\infty}\lambda_{Xm}^{-1}\langle
g,\theta_{m}\varphi_{j}\rangle\theta_{m}\varphi_{j}, g\in
G_{XX}\Biggr\},
\]
we find that $\Gamma_{XX}$ is a one-to-one mapping from the vector
space $%
G_{XX}^{-1}\subset L_{2}(T_{1}\times T_{2})$ onto the vector space
$G_{XX}.$ Thus, restricting $\Gamma_{XX}$ to a subdomain defined by
the subspace $G_{XX}^{-1},$ we can define its inverse for $g\in
G_{XX}$ as
%
%
\begin{equation}
\Gamma_{XX}^{-1}g=\sum_{m,j=1}^{\infty}\lambda
_{Xm}^{-1}\langle g,\theta_{m}\varphi_{j}\rangle\theta_{m}\varphi
_{j}.
\end{equation}
$\Gamma_{XX}^{-1}$ then satisfies the usual properties of an
inverse, in the sense that $\Gamma_{XX}\Gamma_{XX}^{-1}g=g$ for
all $g\in G_{XX},$ and $\Gamma_{XX}^{-1}\Gamma_{XX}h=h$ for
all $h\in G_{XX}^{-1}.$

The following Condition \ref{condc1} for processes $(X,Y)$ is of interest.

\renewcommand{\thecondition}{(C\arabic{condition})}
\setcounter{condition}{0}
\begin{condition}\label{condc1}
The $L_{2}$ -processes $X$%
and $Y$ with Karhunen--Lo\`{e}ve decompositions
(\ref{kl}) satisfy
\[
\sum_{m,j=1}^{\infty}\biggl\{ \frac{E[\xi_{m}\zeta
_{j}]}{\lambda_{Xm}}\biggr\}^2 <\infty.
\]
\end{condition}

 If \ref{condc1} is satisfied, then the solution to the non-invertibility problem
as outlined above is viable in the functional case, as demonstrated
by the following basic result on functional linear models.

\setcounter{thm}{2}
\begin{thm}[(Basic theorem for functional linear models)]\label{th2.3}
A unique solution of the linear model (\ref{linear2}) exists
in $\operatorname{ker}(\Gamma_{XX})^{\perp}$ if and only if
$\mathit{X}$ and $\mathit{Y}$ satisfy Condition
\textup{\ref{condc1}}. In this case, the unique solution is of the form
%
\begin{equation}
\label{regsol}
\beta^\ast_{0}(t,s)=(\Gamma_{XX}^{-1}r_{XY})(t,s).
\end{equation}
\end{thm}

As a consequence of Proposition \ref{pr2.2}, solutions of the functional
linear model (\ref{linear2}), solutions of the functional population
normal equation (\ref{norm}) and minimizers of
$E\|Y-\mathcal{L}_X\beta\|^{2}$ are all equivalent and allow the
usual projection interpretation.

\begin{propositionn}\label{pr2.4}
Assume $X$ and $Y$ satisfy
Condition \textup{\ref{condc1}}. The following are then equivalent:
\begin{enumerate}
\item[(a)] the set of all solutions of the functional linear model
(\ref{linear2});
\item[(b)] the set of all solutions of the population normal equation
(\ref{norm});
\item[(c)] the set of all minimizers of $E\|Y-\mathcal{L}_X\beta
\|^{2}$ for $\beta\in L_{2}(T_{1}\times T_{2})$;
\item[(d)] the set ${\beta_{0}^\ast}+\operatorname{ker}(\Gamma
_{XX})=\{{\beta
_{0}^\ast}+ h | h \in L_{2}(T_{1}\times T_{2}),\Gamma_{XX}h=0\}$.
\end{enumerate}
\end{propositionn}

It is well known that in a finite-dimensional situation, the linear
model (\ref{norm}) always has a unique solution in the column space
of $ \Gamma_{XX}$, which may be obtained by using a generalized
inverse of the matrix $\Gamma_{XX}$. However, in the
infinite-dimensional case, such a solution does not always exist. The
following example demonstrates that a pair of $L_{2}$-processes
does not necessarily satisfy Condition \ref{condc1}. In this case, the
linear model (\ref{norm}) does not have a solution.
\begin{example}\label{ex2.5}
Assume processes $X$ and $Y$ have
Karhunen--Lo\`{e}ve expansions (\ref{kl}), where the random variables
$\xi_{m}$, $\zeta_{j}$ satisfy
%
%
\begin{equation}\label{ex1}
\lambda_{Xm}=E[\xi_{m}^{2}]=\frac{1}{m^{2}}, \qquad \lambda
_{Yj}=E[\zeta_{j}^{2}]=\frac{1}{j^{2}}
\end{equation}
and let
%
%
\begin{equation}\label{ex2}
E[\xi_{m}\zeta_{j}]=\frac{1}{(m+1)^{2}(j+1)^{2}} \qquad \mathrm{for}
\ m,j\geq1.
\end{equation}
As shown in He \textit{et al.}~(\citeyear{He03}), (\ref{ex1}) and (\ref{ex2})
can be
satisfied by a pair of $L_{2}$-processes with appropriate operators
$R_{XX}$, $R_{YY}$ and $R_{XY}$. Then
\begin{eqnarray*}
\sum_{m,j=1}^{\infty}\biggl\{ \frac{E[\xi_{m}\zeta
_{j}]}{\lambda_{Xm}}\biggr\} ^{2}&=&\lim_{n\rightarrow\infty
}\sum_{m,j=1}^{n}\biggl[
\frac{m}{(m+1)(j+1)}\biggr] ^{4}\\
&=&\lim_{n\rightarrow\infty
}\sum_{m=1}^{n}\biggl[ \frac{m}{(m+1)}\biggr]
^{4}\sum_{j=1}^{\infty}\frac{1}{(j+1)^{4}}=\infty
\end{eqnarray*}
and, therefore, Condition \ref{condc1} is not satisfied.
\end{example}

\section{Canonical regression analysis}\label{sec3}

Canonical analysis is a time-honored tool for studying the
dependency between the components of a pair of random vectors or
stochastic processes; for multivariate stationary time series, its
utility was established in the work of Brillinger (\citeyear{Brillinger85}). In
this section, we demonstrate that functional canonical decomposition
provides a useful tool to represent functional linear models. The
definition of functional canonical correlation for $L_{2}$-processes
is as follows.

\begin{definition}\label{def3.1}
The first canonical correlation
$\rho_{1}$ and weight functions $u_{1}$ and $v_{1}$ for $L_{2}$-processes
$X$ and $Y$ are defined as
%
\begin{equation}\label{def1}
\rho_{1}=\sup_{u\in L_{2}(T_{1}), v\in
L_{2}(T_{2})}\operatorname{cov}(\langle u,X\rangle,\langle
v,Y\rangle)=\operatorname{cov}(\langle u_{1},X\rangle,\langle
v_{1},Y\rangle),
\end{equation}
where $u$ and $v$ are subject to
%
\begin{equation}\label{def2} \operatorname{var}(\langle
u_{j},X\rangle)=1,\qquad \operatorname{var}(\langle v_{j},Y\rangle)=1
\end{equation}
for $j=1$. The $k$th canonical correlation $\rho_{k}$ and
weight functions $u_{k}$, $%
v_{k}$ for processes $X$ and $Y$ for $k>1$ are defined as
\[
\rho_{k}=\sup_{u\in L_{2}(T_{1}), v\in L_{2}(T_{2})}
\operatorname{cov}(\langle u,X\rangle,\langle v,Y\rangle
)=\operatorname{cov}(\langle u_{k},X\rangle,\langle
v_{k},Y\rangle),
\]
where $u$ and $v$ are subject to (\ref{def2}) for $j=k$ and
\[
\operatorname{cov}(\langle u_{k},X\rangle,\langle u_{j},X\rangle
)=0,\qquad \operatorname{cov}(\langle
v_{k},Y\rangle,\langle v_{j},Y\rangle)=0
\]
for $j=1,\ldots,k-1.$ We refer to $U_{k}=\langle
u_{k},X\rangle$
and $V_{k}=\langle v_{k},Y\rangle$ as the $k$th canonical variates and
to $(\rho_{k},u_{k},v_{k},U_{k},V_{k})$ as the $k$th canonical
components.
\end{definition}

It has been shown in He \textit{et al.}~(\citeyear{He03}) that canonical
correlations do
not exist for all $L_{2}$-processes, but that Condition \ref{condc2} below
is sufficient for the existence of canonical correlations and weight
functions. We remark that Condition \ref{condc2} implies Condition \ref{condc1}.

\begin{condition}\label{condc2} Let $X$ and $Y$ be
$L_{2}$-processes, with Karhunen--Lo\`{e}ve decompositions
(\ref{kl}) satisfying
\[
\sum_{m,j=1}^{\infty}\biggl\{ \frac{E[\xi_{m}\zeta
_{j}]}{\lambda_{Xm}\lambda_{Yj}^{1/2}}\biggr\} ^{2}<\infty.
\]
\end{condition}

The proposed functional canonical regression analysis exploits
features of functional principal components and of functional
canonical analysis. In functional principal component analysis, one
studies the structure of an $L_{2}$-process via its decomposition
into the eigenfunctions of its autocovariance operator, the
Karhunen--Lo\`{e}ve decomposition (Rice and Silverman (\citeyear{Rice91})). In
functional canonical analysis, the relation between a pair of
$L_{2}$-processes is analyzed by decomposing the processes into
their canonical components. The idea of canonical regression
analysis is to expand the regression parameter function in terms of
functional canonical components for predictor and response
processes. The canonical regression decomposition (Theorem \ref{th3.3})
below provides insights into the structure of the regression
parameter functions and not only aids in the understanding of
functional linear models, but also leads to promising estimation
procedures for functional regression analysis. The details of these
estimation procedures will be discussed in Section \ref{sec4}. We demonstrate
in Section \ref{sec5} that these estimates can lead to competitive prediction
errors in a finite-sample situation.

We now state two key results. The first of these (Theorem \ref{th3.2})
provides the canonical decomposition of the cross-covariance
function of processes $X$ and $Y$. This result plays a central role
in the solution of the population normal equation (\ref{norm}). This
solution is referred to as \textit{canonical regression decomposition}
and it
leads to an explicit representation of the underlying regression
parameter function $\beta_{0}^{\ast}(\cdot,\cdot)$ of the
functional linear model (\ref{linear2}). The decomposition is in
terms of functional canonical correlations $\rho_{j}$ and canonical
weight functions $u_{j}$ and $v_{j}$. Given a predictor process
$X(t)$, we obtain, as a consequence, an explicit representation for
$E(Y(t)|X)=(\mathcal{L}_X\beta_{0}^{\ast})(t)$, where
$\mathcal{L}_X$ is as in (\ref{linear2}). For the following main
results, we refer to the definitions of $\rho_{j}$, $u_{j}$,
$v_{j}$, $U_{j}$, $V_{j}$ in Definition \ref{def3.1}. All proofs are found in
Section \ref{sec7}.

\setcounter{thm}{1}
\begin{thm}[(Canonical decomposition of
cross-covariance function)]\label{th3.2}
Assume that $L_{2}$-processes $X$ and $Y$
satisfy Condition \textup{\ref{condc2}}. The cross-covariance function $r_{XY}$
then allows the following representation in terms of canonical
correlations $\rho_{j}$ and weight functions $u_{j}$ and $v_{j}$:
%
%
\begin{equation}\label{can1}
r_{XY} (s,t)=\sum_{m=1}^{\infty}\rho
_{m}R_{XX}u_{m}(s)R_{YY}v_{m}(t).
\end{equation}
\end{thm}

\begin{thm}[(Canonical regression decomposition)]\label{th3.3}
Assume that the $L_{2}$-processes $X$ and $Y$ satisfy Condition
\textup{\ref{condc2}}. One then obtains, for the regression parameter function $\beta
_{0}^{\ast}(\cdot,\cdot)$ (\ref{regsol}), the following explicit
solution:
%
%
\begin{equation}\label{canreg1}
\beta_{0}^{\ast}(s,t)=\sum_{m=1}^{\infty}\rho
_{m}u_{m}(s)R_{YY}v_{m}(t).
\end{equation}
\end{thm}

 To obtain the predicted value of the response process $Y$, we use
the linear predictor
%
%
\begin{equation}\label{canreg2}
 Y^{\ast}(t)=E(Y(t)|X)=(\mathcal{L}_X\beta_{0}^{\ast
})(t)=\sum_{m=1}^{\infty}\rho_{m}U_{m}R_{YY}v_{m}(t).
\end{equation}

This canonical regression decomposition leads to approximations of
the regression parameter function $\beta_{0}^{\ast}$ and the
predicted process $Y^{\ast}(t)=\mathcal{L}_X\beta_{0}^{\ast}$ via
a finitely truncated version of the canonical expansions
(\ref{canreg1}) and (\ref{canreg2}). The following result provides
approximation errors incurred from finite truncation. Thus, we have a
vehicle to achieve practically feasible estimation of $\beta
_{0}^{\ast}$ and associated predictions $Y^{\ast}$ (Section \ref{sec4}).

\begin{thm}\label{th3.4}
For $K\geq1$, let $\beta
_{K}^{\ast}(s,t)=\sum_{k=1}^{K}\rho_{k}u_{k}(s)R_{YY}v_{k}(t)$ be the
finitely truncated version of the canonical regression decomposition
(\ref{canreg1}) for $\beta_{0}^{\ast}$ and define
$Y_{K}^{\ast}(t)=(\mathcal{L}_X\beta_{K}^{\ast})(t)$. Then,
%
%
\begin{equation} \label{k-pred} Y_{K}^{\ast
}(t)=\sum_{k=1}^{K}\rho_{k}U_{k}R_{YY}v_{k}(t)
\end{equation}
with $E[Y_{K}^{\ast}]=0$. Moreover,
\[
E\|Y^{\ast}-Y_{K}^{\ast}\|^{2}=\sum_{k=K+1}^{\infty}\rho
_{k}^{2}\|R_{YY}v_{k}\|^{2}\rightarrow0 \qquad \mbox{as }  K\rightarrow\infty
\]
and
%
%
\begin{equation}\label{conv} E\|Y-Y_{K}^{\ast
}\|^{2}=E\|Y\|^{2}-E\|\mathcal{L}_X\beta_{K}^{\ast
}\|^{2}=\operatorname{trace}(R_{YY})-\sum_{k=1}^{K}\rho
_{k}^{2}\|R_{YY}v_{k}\|^{2}.
\end{equation}
\end{thm}

In finite-sample implementations, to be explored in the next two
sections, truncation as in (\ref{k-pred}) is a practical necessity;
this requires a choice of suitable truncation parameters.

\section{Estimation procedures}\label{sec4}

\subsection{Preliminaries}\label{sec4.1}

 Estimating the regression parameter function and obtaining
fitted processes from the linear model~(\ref{linear}) based on a
sample of curves is central to the implementation of functional
linear models. In practice, data are observed at discrete time
points and we temporarily assume, for simplicity, that the $N_{x}$
time points are the same for all observed predictor curves and are
equidistantly spaced over the domain of the data. Analogous
assumptions are made for the $N_{y}$ time points where the response
curves are sampled. Thus, the original observations are
$(X_{i},Y_{i})$, $i=1,\ldots,n$, where $X_{i}$ is an
$N_{x}$-dimensional vector sampled at time points $s_j$, and $Y_{i}$ is an
$N_{y}$-dimensional vector sampled at time points $t_j$. We assume
that $N_x$ and $N_y$ are both large. Without going into any
analytical details, we compare the finite-sample behavior of two
functional regression methods, one of which utilizes the canonical
decomposition for regression and the other a well established direct
principal component approach to implement functional linear
regression.

The proposed practical version of functional regression analysis
through functional canonical regression analysis (FCR) is discussed
in Section \ref{sec4.2}. This method is compared with a more standard
functional linear regression implementation that is based on
principal components and referred to as \textit{functional principal
regression} (FPR), in Section \ref{sec4.3}. For the choice of the smoothing
parameters for the various smoothing steps, we adopt
leave-one-curve-out cross-validation (Rice and Silverman (\citeyear{Rice91})).
Smoothing is implemented by local linear fitting for functions and
surfaces (Fan and Gijbels (\citeyear{Fan96})), minimizing locally weighted least
squares.

In a pre-processing step, all observed process data are centered by
subtracting the cross-sectional means $X_{i}(s_{j})-\frac{1}{n}\sum
_{l=1}^{n}X_{l}(s_{j})$, and analogously for $Y_{i}$. If the
data are not sampled on the same grid for different individuals, a
smoothing step may be added before the cross-sectional average is
obtained. As in the previous sections, we use in the following the
notation $X, Y, X_i, Y_i$ to denote centered processes and
trajectories.

When employing the Karhunen--Lo\`{e}ve decomposition (\ref{kl}), we
approximate observed centered processes by the fitted versions
%
\begin{equation}\label{klest} \hat{X}_i(s) = \sum_{l=1}^L \hat\xi_{il}
\hat\theta_l(s),\qquad \hat{Y}_i(t) = \sum_{l=1}^L \hat\zeta_{il}
\hat\varphi_l(t),
\end{equation}
where $\{\hat\theta_l(s)\}_{l=1}^L$ and
$\{\hat\varphi_l(t)\}_{l=1}^L$ are the estimated first $L$ smoothed
eigenfunctions for the random processes $X$ and $Y$, respectively,
with the corresponding estimated eigenscores
$\{\hat\xi_{il}\}_{l=1}^L$ and $\{\hat\zeta_{il}\}_{l=1}^L$ for the
$i$th subject. We obtain these estimates as described in Yao \textit{et
al.}~(\citeyear{Yao05a}). Related estimation approaches, such as those of Rice and
Silverman (\citeyear{Rice91}) or Ramsay and Silverman (\citeyear{Ramsay05}), could alternatively
be used.

\subsection{Functional canonical regression (FCR)}\label{sec4.2}

To obtain functional canonical correlations and the
corresponding weight functions as needed for FCR, we adopt one of the
methods proposed in He \textit{et al.}~(\citeyear{He04}). In preliminary studies, we
determined that the eigenbase method as described there yielded the
best performance for regression applications, with the Fourier base
method a close second. Adopting the eigenbase method, the
implementation of FCR is as follows:
\begin{longlist}
\item[(i)] Starting with the eigenscore estimates as in (\ref{klest}),
estimated raw functional canonical correlations $\tilde{\rho}_{l}$
and $L$-dimensional weight vectors $\mathbf{\hat{u}} _{l}$,
$\mathbf{\hat{v}}_{l}, l=1,\ldots, L$, are obtained by applying
conventional numerical procedures of multivariate canonical analysis
to the estimated eigenscore vectors $(\hat\xi_{i1},\ldots,
\hat\xi_{iL})'$ and $(\hat\zeta_{i1},\ldots, \hat\zeta_{iL})'$. This
works empirically well for moderately sized values of $L$, as
typically obtained from automatic selectors.

\item[(ii)] Smooth weight function estimates $\tilde{u}_l(t),
\tilde{v}_{l}(t)$ are then obtained as
\[
\tilde{u}_l(t)=\mathbf{\hat{u}} _{l}\bolds{\hat\theta}(t),\qquad
\tilde{v}_l(t)=\mathbf{\hat{v}} _{l}\bolds{\hat\varphi}(t),
\]
where
$\bolds{\hat\theta}(t) = (\hat\theta_1(t), \ldots,
\hat\theta_{L}(t))', \bolds{\hat\varphi}(t) = (\hat\varphi_1(t),
\ldots, \hat\varphi_{L}(t))'.$

\item[(iii)] The estimated regression parameter function
$\hat{\beta}$ is obtained according to (\ref{canreg1}) by
\[
\hat{\beta
}(s,t)=\sum_{l=1}^{L}\tilde{\rho}_{l}\mathbf{\tilde{u}}_{l}(s)
\int_{T_2} \hat{r}_{YY}(s,t) \tilde{v}_l(s) \,\mathrm{d}s,
\]
where
$\hat{r}_{YY}$ is an estimate of the covariance function of $Y$,
obtained by two-dimensional smoothing of the empirical
autocovariances of $Y$. This estimate is obtained as described in
Yao \textit{et al.}~(\citeyear{Yao05a}). Since the data are regularly sampled, the above
integrals are easily obtained by the approximations
$\sum_{j=1}^{m_y} \hat{r}_{YY}(t_j,t)
\tilde{v}_l(t_j)(t_{j}-t_{j-1}), l=1,\ldots, L$, with $t_0$
defined analogously to $s_0$ in (\ref{int}) below.

\item[(iv)] Fitted/predicted processes
%
%
\begin{equation}\label{pred}
\hat{Y}_{i}(t)=\int_{T_{1}}\hat{\beta}(s,t)X_{i}(s)\,\mathrm{d}s
\qquad\mbox{for } i=1,\ldots,n,
\end{equation}
are obtained, where the integral is again evaluated numerically by
%
\begin{equation}\label{int}
\hat{Y}_{i}(t)=\sum_{j=1}^{N_{x}}\hat{\beta}
(s_{j},t)X_{i}(s_{j})(s_{j}-s_{j-1}).
\end{equation}
Here, $s_{0}$ is chosen such that $s_{1}-s_{0}=s_{2}-s_{1}$.
\end{longlist}

This procedure depends on two tuning parameters, a bandwidth $%
h $ for the smoothing steps (which are defined in detail, e.g., in Yao
\textit{et al.}~(\citeyear{Yao05a})) and the number of canonical \mbox{components}~$L$ that are
included. These tuning parameters may be determined by leave-one-out
cross-validation (Rice and Silverman (\citeyear{Rice91})) as follows. With $\alpha
=(h,L)$, the $i$th leave-one-out estimate for $\beta$ is
%
%
\begin{equation}\label{cvbet}
\hat{\beta}_{\alpha}^{(-i)}=\sum_{l=1}^{L}\breve{\rho}
_{l}^{(-i)}\tilde{u}_{h,l}^{(-i)}(s)\hat{R}_{YY}^{(-i)}(\tilde{v}
)_{h,l}^{(-i)}(t)\qquad \mbox{for } i=1,\ldots,n,
\end{equation}
where $\breve{\rho}_{l}^{(-i)}$ is the $l$th canonical correlation,
and $ \tilde{u}_{h,l}^{(-i)}$ and
$\hat{R}_{YY}^{(-i)}\tilde{v}_{h,l}^{(-i)}$ are the $l$th weight
function untransformed and transformed with the covariance operator,
respectively, all obtained while leaving out the data for the $i$th
subject. Computation of these estimates
follows steps (iii) and (iv) above, using tuning parameter $\alpha
=(h,L)$, and omitting the $i$th pair of observed curves
$(X_{i},Y_{i})$. The average leave-one-out squared prediction error
is then
%
%
\begin{equation}\label{pe}
\mathit{PE}_{\alpha}=\frac{1}{n}\sum_{i=1}^{n}\int_{T_2}
\biggl(Y_{i}(t)-\int_{T_1} X_{i}(s) \hat{\beta}_{\alpha
}^{(-i)}(s,t)\,\mathrm{d}s\biggr)^{2}\,\mathrm{d}t.
\end{equation}
The cross-validation procedure then selects the tuning parameter
that minimizes the approximate average prediction error,\vspace*{-2pt}
\begin{eqnarray*}
\hat{\alpha}=\mathop{\arg\min}_{\alpha}\tilde{\mathit{PE}}_{\alpha},
\end{eqnarray*}
where $\tilde{\mathit{PE}}_{\alpha}$ is obtained by replacing the
integrals on the right-hand side of (\ref{pe}) by sums of the type
(\ref{int}).

\subsection{Functional principal component regression (FPR)}\label{sec4.3}

Yao \textit{et al.}~(\citeyear{Yao05b}) considered an implementation of functional linear
regression whereby one uses functional principal component analysis
for predictor and response functions separately, followed by simple
linear regressions of the response principal component scores on the
predictor scores. We adopt this approach as FPR.

Briefly, defining $\sigma_{mp}=E(\xi_m\zeta_p)$, this approach is
based on representations\vspace*{-2pt}
\[
\beta(s,t) = \sum_{m=1}^\infty\sum_{p=1}^\infty
\frac{\sigma_{mp}}{{\lambda}_{Xm}}\theta_m(s)\varphi_p(t)
\]
of the
regression parameter function $\beta(s,t)$, where\vspace*{-2pt}
%
\begin{eqnarray}\label{sig}
\sigma_{mp} = \int_{T_2}\int_{T_1} \theta_m(s) r_{XY}(s,t)
\varphi_p(t) \,\mathrm{d}s \,\mathrm{d}t
\end{eqnarray}
for all $m$ and $p$.

For estimation, one first obtains a smooth estimate $\hat{r}_{XY}$
of the cross-covariance $r_{XY}$ by smoothing sample
cross-covariances, for example, by the method described in Yao \textit
{et al.}~(\citeyear{Yao05b}). This leads to estimates $\hat{\sigma}_{mp}$ of
$\sigma_{mp}, 1 \le m,p \le L,$ by plugging in estimates
$\hat{r}_{XY}$ for $r_{XY}$ and $\hat{\theta}_l, \hat{\varphi}_l$
for eigenfunctions $\theta_l, \varphi_l$ (as described in Section \ref{sec4.2}),
in combination with approximating the integrals in (\ref{sig})
by appropriate sums. One may then use these estimates in conjunction
with estimates~$\hat{\lambda}_{Xm}$ of eigenvalues $\lambda_{Xm}$
to arrive at the estimate $\hat\beta$ of the regression parameter
function~$\beta(s,t)$ given by\vspace*{-2pt}
\[
\hat\beta(s,t) = \sum_{m=1}^L \sum_{p=1}^L
\frac{\hat\sigma_{mp}}{\hat{\lambda}_{Xm}}\hat\theta_m(s)\hat
\varphi_p(t).
\]
For further details about numerical implementations, we refer to Yao
\textit{et al.}~(\citeyear{Yao05b}).\vadjust{\goodbreak}

\section{Application to medfly mortality data}\label{sec5}

In this section, we present an application to age-at-death data
that were collected for cohorts of male and female medflies in a
biodemographic study of survival and mortality patterns of cohorts
of male and female Mediterranean fruit flies (\textit{Ceratitis
capitata}; for details, see Carey \textit{et al.}~(\citeyear{Carey02})). A point of interest
in this study is the relation of mortality trajectories between male
and female medflies which were raised in the same cage. One
specifically desires to quantify the influence of male survival on female
survival. This is of interest because female survival determines the
number of eggs laid and thus reproductive success of these flies. We
use a subsample of the data generated by this experiment, comprising
46 cages of medflies, to address these questions. Each cage contains
both a male and a female cohort, consisting each of approximately
4000 male and 4000 female medflies. These flies were raised in the
shared cage from the time of eclosion. For each cohort, the number
of flies alive at the beginning of each day was recorded, simply by
counting the dead flies on each day; we confined the analysis to the
first 40 days. The observed processes $X_{i}(t)$ and $Y_{i}(t)$,
$t=1,\ldots,40, i=1,\ldots, 46,$ are the estimated random hazard
functions for male and female cohorts, respectively. All deaths are
fully observed so that censoring is not an issue. In a
pre-processing step, cohort-specific hazard functions were estimated
nonparametrically from the lifetable data, implementing the
transformation approach described in M\"{u}ller \textit{et al.}~(\citeyear{Muller97a}).

\begin{table}[b]
\caption{Results for medfly data, comparing functional
canonical regression (FCR) and functional principal component
regression (FPR) with regard to average  leave-one-out squared
prediction error (PE) (\protect\ref{pe}); values for bandwidth $h$ and
number $L$ of components as chosen by cross-validation are also
shown}\label{t1}
\begin{tabular}{@{}llll@{}}
\hline
&  \multicolumn{1}{l}{$h$}  &  \multicolumn{1}{l}{$L$} &  \multicolumn{1}{l}{PE}\\
\hline
FCR &  1.92  & 3 & 0.0100 \\
FPR & 1.65 & 3 &  0.0121  \\
\hline
\end{tabular}
\end{table}

A functional linear model was used to study the specific influence
of male mortality on female mortality for flies that were raised in the same
cage, with the hazard function of males as predictor process and
that of females as response process. We applied both the proposed
regression via canonical representation (FCR) and the more
conventional functional regression based on principal components
(FPR), implementing the estimation procedures described in the
previous section. Tuning parameters were selected by
cross-validation. Table \ref{t1} lists the average squared prediction
error (PE) (\ref{pe}) obtained by the leave-one-out technique. For
this application, the FCR procedure is seen to perform about 20\%
better than FPR in terms of PE.

\begin{figure}

\includegraphics{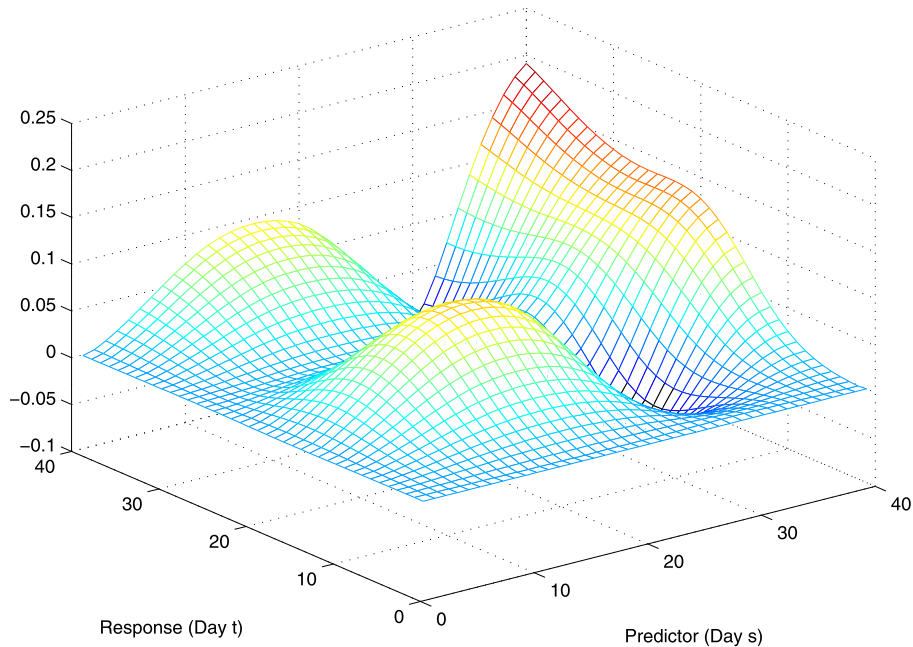}
\caption{Estimated regression parameter surface
obtained by functional canonical regression for the medfly study.}\label{f1}
\end{figure}

The estimated regression parameter surface $\hat{\beta}(s,t)$
that is obtained for the FCR regression when choosing the
cross-validated values for $h$ and $L$, as given in Table \ref{t1}, is
shown in Figure \ref{f1}. The shape of the regression surface indicates
that female mortality at later ages is very clearly affected by male
mortality throughout male lifespan, while female mortality at very
early ages is not much influenced by male mortality. The effect of
male mortality on female mortality is periodically elevated, as
evidenced by
the bumps visible in the surface. The particularly influential
predictive periods are male mortality around days 10 and 20, which
then has a particularly large influence on female mortality around
days 15 and 25, that is, about five days later, and, again, around days 35
and 40, judging from the locations of the peaks in the surface of
$\hat{\beta}(s,t)$. In contrast, enhanced male mortality around
day 30 leads to lessened female mortality throughout, while enhanced
male mortality at age 40 is associated with higher older-age female
mortality. These observations point to the existence of periodic
waves of mortality, first affecting males and subsequently females.
While some of the waves of increased male mortality tend to be
associated with subsequently increased female mortality, others are
associated with subsequently decreased female mortality.

These waves of mortality might be related to the so-called
``vulnerable periods'' that are characterized by locally heightened
mortality (M\"{u}ller \textit{et al.}~(\citeyear{Muller97b})). One such vulnerable
period occurs
around ages 10 and 20, and the analysis suggests that heightened
male mortality during these phases is indicative of heightened
female mortality. In contrast, heightened male mortality during a
non-vulnerable period such as the time around 30 days seems to be
associated with lower female mortality. A word of caution is in
order as no inference methods are available to establish that the
bumps observed in $\hat{\beta}(s,t)$ are real, so one cannot
exclude the possibility that these bumps are enhanced by random
fluctuations in the data.

\begin{figure}

\includegraphics{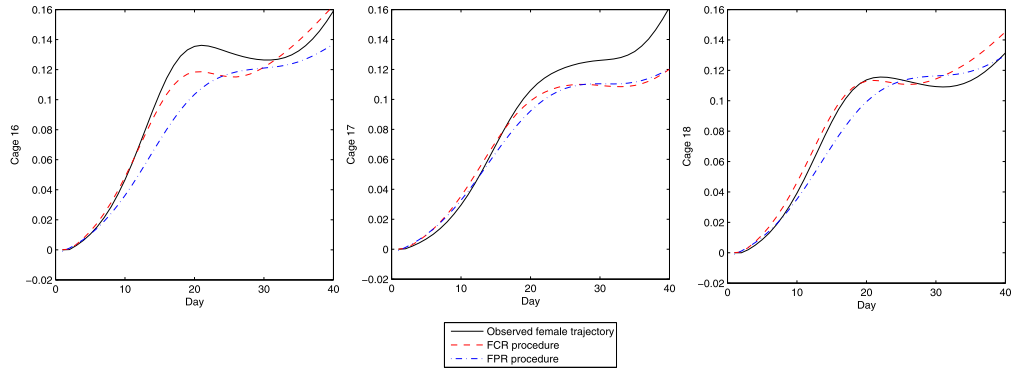}
\caption{Functional regression of female (response)
on male (predictor) medfly trajectories quantifying  mortality in
the form of cohort hazard functions for three cages of flies. Shown
are actually observed female trajectories (solid) that are not used
in the prediction, as well as the predictions for these trajectories
obtained through estimation procedures based on functional principal
component regression (FPR) (dash--dot) and on functional canonical
regression (FCR) (dashed).}\label{f2}
\end{figure}

Examples of observed, as well as predicted, female mortality
trajectories for three randomly selected pairs of cohorts (male and
female flies raised in the same cages) are displayed in Figure~\ref{f2}.
The predicted female trajectories were constructed by applying both
regression methods (FCR and FPR) with the leave-one-out technique.
The prediction of an individual response trajectory from a predictor
trajectory cannot, of course, be expected to be very close to the
actually observed response trajectory, due to the extra random
variation that is a large inherent component of response
variability; this is analogous to the situation of predicting an
individual response in the well-known simple linear regression
case. Nevertheless, overall, FCR predictions are found to be closer
to the target.

We note the presence of a ``shoulder'' at around day 20 for the
three female mortality curves. This ``shoulder'' is related to the
wave phenomenon visible in $\hat{\beta}(s,t)$ as discussed
above and corresponds to a phase of elevated female mortality. The
functional regression method based on FCR correctly predicts the
shoulder effect and its overall shape in female mortality. At the
rightmost points, for ages near 40 days, the variability of the
mortality trajectories becomes large, posing extra difficulties for
prediction in the right tail of the trajectories.

\section{Additional results}\label{sec6}

Theorems \ref{th6.3} and
\ref{th6.4} in this section provide a functional analog to the sums-of-squares
decomposition of classical regression analysis. In addition,
we provide two results characterizing the regression operators
$\mathcal{L}_X$. We begin with two auxiliary results which are taken
from He \textit{et al.}~(\citeyear{He03}). The first of these characterizes the
correlation operator between processes $X$ and $Y$.

\begin{lemma}\label{le6.1} Assume that the
$L_{2}$-processes $X$ and
$Y$ satisfy Condition \textup{\ref{condc2}}. The correlation operator
$R_{XX}^{-1/2}R_{XY}R_{YY}^{-1/2}$can then be extended continuously to a
Hilbert--Schmidt operator $R$ on $L_{2}(T_{2})$ to $L_{2}(T_{1})$.
Hence, $R_{0}=R^{\ast}R$ is also a Hilbert--Schmidt operator with a
countable number of non-zero eigenvalues and
eigenfunctions $\{(\lambda_{m},q_{m})\}$, $m\geq1, \lambda_{1}\geq
\lambda_{2}\geq\cdots$,  $%
p_{m}=Rq_{m}/\sqrt{\lambda_{m}}.$ Then:
\begin{longlist}
\item[(a)] $\rho_{m}=\sqrt{\lambda_{m}}$,
$u_{m}=R_{XX}^{-1/2}p_{m},$ $%
v_{m}=R_{YY}^{-1/2}q_{m}$ and both $u_{m}$ and $v_{m}$ are
$L_{2}$-functions;
\item[(b)] $\operatorname{corr}(U_{m},U_{j})=\langle
u_{m},R_{XX}u_{j}\rangle=\langle
p_{m},p_{j}\rangle=\delta_{mj}$;
\item[(c)]$\operatorname{corr}(V_{m},V_{j})=\langle
v_{m},R_{XX}v_{j}\rangle=\langle
q_{m},q_{j}\rangle=\delta_{mj}$;
\item[(d)]$\operatorname{corr}(U_{m},V_{j})=\langle
u_{m},R_{XX}v_{j}\rangle=\langle
p_{m},Rq_{j}\rangle=\rho_{m}\delta_{mj}$.
\end{longlist}
\end{lemma}

One of the main results in He \textit{et al.}~(\citeyear{He03}) reveals that the
$L_{2}$%
-processes $X$ and $Y$ can be expressed as sums of uncorrelated
component functions and the correlation between the $m$th components
of the expansion is the $m$th corresponding functional canonical
correlation between the two processes.

\begin{lemma}[(Canonical decomposition)]\label{le6.2}  Assume
$L_{2}$-processes $X$ and $Y$ satisfy Condition \textup{\ref{condc2}}. There then
exists a decomposition:
\begin{longlist}
\item[(a)]
\[
X=X_{c,K}+X_{c,K}^{\bot},\qquad Y=Y_{c,K}+Y_{c,K}^{\bot},
\]
where
\begin{eqnarray*}
X_{c,K}&=&\sum_{j=1}^{K}U_{j}R_{XX}u_{j},\qquad  X_{c,K}^{\bot
}=X-X_{c,K},\\
Y_{c,K}&=&\sum_{j=1}^{K}V_{j}R_{YY}v_{j},\qquad
Y_{c,K}^{\bot}=Y-Y_{c,K}.
\end{eqnarray*}
The index $K$ stands for canonical decomposition with $K$
components,
and $U_{j},$ $V_{j},$ $u_{j},$ $v_{j}$ are as in Definition \ref{def3.1}. Here, $
(X,Y) $ and $(X_{c,K},Y_{c,K})$ share the same first $K$ canonical
components, and $(X_{c,K},Y_{c,K})$ and $(X_{c,K}^{\bot
},Y_{c,K}^{\bot})$ are uncorrelated, that is,
\begin{eqnarray*}
\operatorname{corr}(X_{c,K},X_{c,K}^{\bot})&=&0,\qquad \operatorname
{corr}(Y_{c,K},Y_{c,K}^{\bot
})=0,\\
 \operatorname{corr}(X_{c,K},Y_{c,K}^{\bot})&=&0,\qquad
\operatorname{corr}(Y_{c,K},X_{c,K}^{\bot})=0.
\end{eqnarray*}
\item[(b)] Let $K\rightarrow\infty$ and $X_{c,\infty}=
\sum_{m=1}^{\infty}U_{m}R_{XX}u_{m}, Y_{c,\infty
}=\sum_{m=1}^{\infty}V_{m}R_{YY}v_{m}. $ Then
\[
X=X_{c,\infty}+X_{c,\infty}^{\bot}, \qquad  Y=Y_{c,\infty}+%
X_{c,\infty}^{\bot},
\]
where $X_{c,\infty}^{\bot}=X-$ $X_{c,\infty},$ $Y_{c,\infty
}^{\bot}=Y-$ $Y_{c,\infty}$. Here, $(X_{c,\infty},Y_{c,\infty})$
and $(X,Y)$ share the same canonical components, $\operatorname
{corr}(X_{c,\infty
}^{\bot},Y_{c,\infty}^{\bot})=0$, and $(X_{c,\infty}^{\bot
},Y_{c,\infty}^{\bot})$ and $(X_{c,\infty},Y_{c,\infty})$ are
uncorrelated. Moreover, $X_{c,\infty}^{\bot}=0$ if $\{p_{m},m\geq
1\}$ forms a basis of the closure of the domain of $R_{XX}$ and
$Y_{c,\infty}^{\bot}=0$ if $\{q_{m},m\geq1\}$ forms a basis of
the closure of the domain of $R_{YY}$.
\end{longlist}
\end{lemma}

Since the covariance operators of $L_{2}$-processes are non-negative
self-adjoint, they can be ordered as follows. The definitions of
$Y^{\ast}, Y_{K}^{\ast}, Y_{c,\infty}$ are in (\ref{canreg2}),
(\ref{k-pred}) and Lemma \ref{le6.2}(b), respectively.

\begin{thm}\label{th6.3} For $K\geq1,$
$R_{Y_{K}^{\ast
}Y_{K}^{\ast}}\leq R_{Y^{\ast}Y^{\ast}}\leq R_{Y_{c,\infty
}Y_{c,\infty}}\leq R_{YY}$.
\end{thm}

In multiple regression analysis, the ordering of the operators in
Theorem \ref{th6.3} is related to the ordering of regression models in terms
of a notion analogous to the regression sum of squares (SSR). The
canonical regression decomposition provides information about the
model in terms of its canonical components. Our next result
describes the canonical correlations between observed and fitted
processes. This provides an extension of the coefficient of multiple
determination, $R ^{2}=\operatorname{corr}(Y,\hat{Y}),$ an
important quantity
in classical multiple regression analysis, to the functional case;
compare also Yao \textit{et al.}~(\citeyear{Yao05b}).

\begin{thm}\label{th6.4}
Assume that $L_{2}$-processes
$X$ and
$Y$ satisfy Condition \textup{\ref{condc2}}. The canonical correlations and
weight functions for the pair of observed and fitted response
processes $(Y,Y^{\ast})$ are then $\{(\rho_{m},v_{m},v_{m}/\rho_{m});
m\geq1\}$ and the corresponding $K$%
-component (or $\infty$-component) canonical decomposition for
$Y^{\ast},$
as defined in Lemma \ref{le6.2} for $K \ge1$ and denoted here by
$Y_{c,K}^{\ast
}$ (or %
$Y_{c,\infty}^{\ast}$), is equivalent to the process $Y_{K}^{\ast} $
or $Y^{\ast}$ given in Theorem \ref{th3.4}, that is,
%
%
\begin{equation}\label{t6.4}
Y_{c,K}^{\ast}=Y_{K}^{\ast}=\sum_{m=1}^{K}\rho
_{m}U_{m}R_{YY}v_{m},\qquad K\geq1,\qquad Y_{c,\infty}^{\ast}=Y^{\ast
}=\sum_{m=1}^{\infty}\rho_{m}U_{m}R_{YY}v_{m}.
\end{equation}
\end{thm}

 We note that if $Y$ is a scalar, then $R^{2}=\rho_{1},$ and for a
functional response $Y$, $R^{2}$ is replaced by the set $\{\rho
_{m},$ $m\geq1\}.$

The following two results serve to characterize the
regression operator $\mathcal{L}_X$ defined in (\ref{linear2}). They
are used in the proofs provided in the following section.

\begin{propositionn}\label{pr6.5}The adjoint operator of
$\mathcal{L}_X$ is $\mathcal{L}_X^{\ast}\dvtx L_{2}(T_{2})\rightarrow
L_{2}(T_{1}\times T_{2}) $, where
\[
(\mathcal{L}_X^{\ast}z)(s,t)=X(s)z(t)\qquad   \mbox{for }  z\in L_{2}(T_{2}).
\]
\end{propositionn}

We have the following relation between the correlation operator $
\Gamma_{XX}$ defined in (\ref{gamma}) and the regression operator
$\mathcal{L}_X$.

\begin{propositionn}\label{pr6.6} The operator $\Gamma
_{XX}$ is a
self-adjoint non-negative Hilbert--Schmidt operator and satisfies $
\Gamma_{XX}=E[\mathcal{L}_X^{\ast}\mathcal{L}_X].$
\end{propositionn}

\section{Proofs}\label{sec7}

In this section, we provide sketches of proofs and some auxiliary
results. We use tensor notation to define an operator $\theta
\otimes\varphi\dvtx H\rightarrow H,$
\[
(\theta\otimes\varphi) (h)=\langle h,\theta\rangle\varphi\qquad\mbox{for}
\
h\in H.
\]

\begin{pf*}{Proof of Proposition \protect\ref{pr2.2}} To prove (a) $
\Rightarrow$
(b), we multiply equation (\ref{linear2}) by $X$ on both sides and
take expected values to obtain $E(XY)=E(X\mathcal{L}_{X}\beta
_{0})+E(X\varepsilon)$. Equation (\ref{norm}) then follows from
$E(XY)=r_{XY},$ $E(X\mathcal{L}_{X}\beta
_{0})=\Gamma_{XX}\beta_{0}$ (by Propositions \ref{pr6.5} and \ref{pr6.6}) and $%
E(X\varepsilon)=0$.

For (b) $ \Rightarrow$ (c), let $\beta_{0}$ be a solution of
equation (\ref{norm}). For any $\beta\in L_{2}(T_{1}\times
T_{2})$, we then have
$E\|Y-\mathcal{L}_X\beta\|^{2}=E\|Y-\mathcal{L}_X\beta
_{0}\|^{2}+E\|\mathcal{L}_X(\beta_{0}-\beta)\|^{2}+2E[\langle
Y-\mathcal{L}_X\beta_{0},\mathcal{L}_X(\beta_{0}-\beta)\rangle].$
Since
\begin{eqnarray*}
&&E\langle Y-\mathcal{L}_X\beta_{0},\mathcal
{L}_X(\beta
_{0}-\beta)\rangle\\
&&\quad=E\langle \mathcal{L}_X^{\ast}Y-\mathcal
{L}_X^{\ast
}\mathcal{L}_X\beta_{0},\beta_{0}-\beta\rangle\\
&&\quad=\langle E[\mathcal{L}_X^{\ast}Y]-E[\mathcal{L}_X^{\ast
}\mathcal{L}_X\beta_{0}],\beta_{0}-\beta\rangle=\langle
r_{XY}-\Gamma
_{XX}\beta_{0},\beta_{0}-\beta\rangle=0,
\end{eqnarray*}
by Proposition \ref{pr6.6},
we then have
\[
E\|Y-\mathcal{L}_X\beta\|^{2}=E\|Y-\mathcal{L}_X\beta
_{0}\|^{2}+E\|\mathcal{L}_X(\beta_{0}-\beta)\|^{2}\geq
E\|Y-\mathcal{L}_X\beta_{0}\|^{2},
\]
which implies that $\beta_{0}$ is indeed a minimizer of
$E\|Y-\mathcal{L}_X\beta\|^{2}$.

For $\mathrm{(c)}\Rightarrow\mathrm{(a)}$, let
\[
d^{2}=E\|Y-\mathcal{L}_X\beta_{0}\|^{2}=\mathop{\min}_{\beta\in
L_{2}(T_{1}\times T_{2})}E\|Y-\mathcal{L}_X\beta\|^{2}.
\]
Then, for any $\beta\in L_{2}(T_{1}\times T_{2}),$ $a\in\mathbf{R,}$
\begin{eqnarray*}
d^{2}&=&E\|Y-\mathcal{L}_X\beta_{0}\|^{2}\leq
E\|Y-\mathcal{L}_X(\beta_{0}+a\beta)\|^{2}
\\
&=&E\|Y-\mathcal{L}_X\beta_{0}\|^{2}-2E\langle Y-\mathcal{L}_X\beta
_{0},\mathcal{L}_X(a\beta)\rangle+E\|\mathcal{L}_X(a\beta)\|^{2}\\
&=&d^{2}-2a\langle E[X(Y-\mathcal{L}_X\beta_{0})],\beta
\rangle+a^{2}E\|\mathcal{L}_X\beta\|^{2}.
\end{eqnarray*}
Choosing $a=\langle E[X(Y-\mathcal{L}_X\beta_{0})],\beta
\rangle/E\|\mathcal{L}_X\beta\|^{2},$ it follows that $|\langle
E[X(Y-\mathcal{L}_X\beta_{0})],\beta
\rangle|^{2}/\break E\|\mathcal{L}_X\beta\|^{2}\leq0$ and $\langle
E[X(Y-\mathcal{L}_X\beta_{0})],\beta\rangle=0.$ Since $\beta$ is
arbitrary, $E[X(Y-\mathcal{L}_X\beta_{0})]=0$ and therefore $\beta
_{0}$ satisfies the functional linear model (\ref{linear2}).
\end{pf*}

\begin{pf*}{Proof of Theorem \protect\ref{th2.3}} Note, first, that
$r_{XY}(s,t)=\sum_{m,j}E[\xi_{m}\zeta_{j}]\theta_{m}(s)\varphi
_{j}(t).$ Thus, Condition~\ref{condc1} is equivalent to $r_{XY}\in G_{XX}$.
Suppose that a unique solution of (\ref{linear2}) exists in
$\operatorname{ker}(\Gamma
_{XX})^{\bot}.$
This solution is then also a solution of (\ref{norm}), by Proposition
\ref{pr2.2}(b). Therefore, $%
r_{XY}\in G_{XX}$, which implies \ref{condc1}. On the other hand, if \ref{condc1}
holds, then $r_{XY}\in G_{XX},$ which implies that $\Gamma
_{XX}^{-1}r_{XY}=\sum_{m}\lambda_{Xm}^{-1}\langle r_{XY},\theta
_{m}\varphi_{j}\rangle\theta_{m}\varphi_{j}$ is a solution of
(\ref{norm}), is in $\operatorname{ker}(\Gamma_{XX})^{\bot}$ and,
therefore, is
the unique solution in $\operatorname{ker}(\Gamma_{XX})^{\bot}$ and
also the
unique solution of (\ref{linear2}) in $\operatorname{ker}(\Gamma
_{XX})^{\bot}.$
\end{pf*}
\begin{pf*}{Proof of Proposition \protect\ref{pr2.4}}
 The equivalence of (a), (b)
and (c) follows from Proposition \ref{pr2.2} and (d) $ \Rightarrow$ (b)
is a consequence of Theorem \ref{th2.3}. We now prove (b) $ \Rightarrow
$ (d). Let $\beta_{0}$ be a solution of~(\ref{norm}). Proposition
\ref{pr2.2}
and Theorem \ref{th2.3} imply that both $\beta_{0}$ and $\beta_{0}^{\ast}$ minimize $
E\|Y-\mathcal{L}_X\beta\|^{2}$ for $\beta\in L_{2}(T_{1}\times
T_{2}).$ Hence, $E\|Y-\mathcal{L}_X\beta
_{0}\|^{2}=E\|Y-\mathcal{L}_X\beta_{0}^{\ast
}\|^{2}+E\|\mathcal{L}_X(\beta_{0}^{\ast}-\beta_{0})\|
^{2}+2E\langle
Y-\mathcal{L}_X\beta_{0}^{\ast},\mathcal{L}_X(\beta_{0}^{\ast
}-\beta_{0})\rangle,$ which, by Proposition \ref{pr6.6}, implies that
$2E\langle
\mathcal{L}_X^{\ast}(Y-\mathcal{L}_X\beta_{0}^{\ast}),\beta
_{0}^{\ast}-\beta_{0}\rangle=2\langle r_{XY}-\Gamma_{XX}\beta
_{0}^{\ast},\beta_{0}^{\ast}-\beta_{0}\rangle=0.$ Therefore,
$E\|\mathcal{L}_X(\beta_{0}^{\ast}-\beta_{0})\|^{2}=\|\Gamma
_{XX}^{1/2}(\beta_{0}^{\ast}-\beta_{0})\|^{2}=0.$ It follows that
$\beta_{0}^{\ast}-\beta_{0}\in\operatorname{ker}(\Gamma_{XX}),$
or $\beta
_{0}=\beta_{0}^{\ast}+h,$ for an
$h\in\operatorname{ker}(\Gamma_{XX}).$\
\end{pf*}

\begin{pf*}{Proof of Theorem \protect\ref{th3.2}} According to Lemma \ref{le6.2}(b),
Condition \ref{condc2} guarantees the existence of the canonical components
and canonical decomposition of $X$ and $Y$. Moreover,
\begin{eqnarray*}
r_{XY}(s,t)&=&E[X(s)Y(t)] =E\bigl[\bigl(X_{c,\infty}(s)+X_{c,\infty
}^{\bot}(s)\bigr)\bigl(Y_{c,\infty}(t)+Y_{c,\infty}^{\bot}(t)\bigr)\bigr]\\
&=&E[X_{c,\infty}(s)Y_{c,\infty}(t)]
=E\Biggl[\sum_{m=1}^{\infty
}U_{m}R_{XX}u_{m}(s)\sum_{m=1}^{\infty}V_{m}R_{YY}v_{m}(t)\Biggr]
\\&=&\sum_{m,j=1}^{\infty
}E[U_{m}V_{j}]R_{XX}u_{m}(s)R_{YY}v_{m}(t)
=\sum_{m=1}^{\infty}\rho
_{m}R_{XX}u_{m}(s)R_{YY}v_{m}(t).
\end{eqnarray*}
We now show that
the exchange of the expectation with the summation above is valid.
From Lemma \ref{le6.1}(b), for any $K>0$ and the spectral decomposition
$R_{XX}=\sum_{m}\lambda_{Xm}\theta_{m}\otimes\theta_{m}$,
\begin{eqnarray*}
\sum_{m=1}^{K}E\|U_{m}R_{XX}u_{m}\|^{2}&=&\sum
_{m=1}^{K}E[U_{m}^{2}]\|R_{XX}^{1/2}p_{m}\|^{2}=\sum
_{m=1}^{K}\langle p_{m},R_{XX}p_{m}\rangle
\\
&=&\sum_{m=1}^{K} \sum_{j=1}^{\infty}\lambda_{Xj}\langle
p_{m},\theta_{j}\rangle^{2}=\sum_{j=1}^{\infty}\lambda
_{Xj}\Biggl(\sum_{m=1}^{K}\langle p_{m},\theta_{j}\rangle^{2}\Biggr)\\
\\
&\leq& \sum_{j=1}^{\infty}\lambda_{Xj}\|\theta
_{j}\|^{2}=\sum_{j=1}^{\infty}\lambda_{Xj}<\infty,
\end{eqnarray*}
where the inequality follows from the fact that
$\sum_{m=1}^{K}\langle p_{m},\theta_{j}\rangle^{2}$ is the square length
of
the projection of $\theta_{j}$ onto the linear subspace spanned by $%
\{p_{1},\ldots,p_{K}\}$. Similarly, we can show that for any $K>0$,
\[
\sum_{m=1}^{K}E\|V_{m}R_{YY}v_{m}\|^{2}<\sum
_{j=1}^{\infty}\lambda_{Yj}<\infty.
\]
\upqed\end{pf*}

\begin{pf*}{Proof of Theorem \protect\ref{th3.3}} Note that Condition \ref{condc2} implies
Condition \ref{condc1}. Hence, from Theorem~\ref{th2.3}, $\beta_{0}^{\ast}=
\Gamma_{XX}^{-1}r_{XY}$ exists and is unique in $%
\operatorname{ker}(\Gamma_{XX})^{\bot}$. We can show (\ref
{canreg1}) by applying
$\Gamma_{XX}^{-1}$ to both sides of (\ref{norm}), exchanging the
order of summation and integration. To establish (\ref{canreg2}), it
remains to show that
%
%
\begin{equation}\label{p1}
\sum_{m=1}^{\infty}\|\rho
_{m}u_{m}R_{YY}v_{m}\|^{2}<\infty,
\end{equation}
where $u_{m}R_{YY}v_{m}(s,t)=u_{m}(s)R_{YY}v_{m}(t)$ in
$L_{2}(T_{1}\times T_{2}).$
Note that
\[
\rho_{m}u_{m}=\rho
_{m}R_{XX}^{1/2}p_{m}=R_{XX}^{1/2}Rq_{m}=\sum_{j=1}^{\infty
}\frac{1}{ \sqrt{\lambda_{Xj}}}\langle Rq_{m},\theta_{j}\rangle
\theta
_{j},
\]
where the operator $R=R_{XX}^{1/2}R_{XY}R_{YY}^{1/2}$ is defined in
Lemma \ref{le6.1} and can be written as $R=\sum_{k,\ell}r_{k\ell}\varphi
_{k}\otimes\theta_{\ell},$ with $r_{km}=E[\xi_{k}\zeta
_{\ell}]/\sqrt{\lambda_{Xk}\lambda_{Y\ell}},$ using the
Karhunen--Lo\`{e}ve expansion (\ref{kl}). Then,
\[
Rq_{m}=\sum_{k,\ell}r_{k\ell}\langle \varphi_{k},q_{m}\rangle
\theta
_{\ell},\qquad \langle Rq_{m},\theta_{j}\rangle=\sum_{k}r_{kj}\langle
\varphi
_{k},q_{m}\rangle
\]
and, therefore,
\begin{eqnarray*}
&&\sum_{m}\|\rho_{m}u_{m}R_{YY}v_{m}\|^{2}\\
&&\quad\leq
\sum_{m}\|\rho_{m}u_{m}\|^{2}\|R_{YY}v_{m}\|^{2}
=\sum_{m}\biggl[\sum_{j}\frac{1}{\lambda_{Xj}}(\langle
Rq_{m},\theta_{j}\rangle)^{2}\biggr] \|R_{YY}v_{m}\|^{2}\\
&&\quad=\sum_{m}\biggl[\sum_{j}\frac{1}{\lambda_{Xj}}\biggl\{
\sum_{k}r_{kj}\langle \varphi_{k},q_{m}\rangle\biggr\}
^{2}\biggr]\|R_{YY}v_{m}\|^{2}\\
&&\quad \leq
\sum_{m}\biggl[\sum_{j}\frac{1}{\lambda_{Xj}}
\sum_{k}r_{kj}^{2}\sum_{\ell}\langle \varphi
_{\ell},q_{m}\rangle^{2}\biggr]\|R_{YY}v_{m}\|^{2}\\
&&\quad=\biggl[\sum_{j}\frac{1}{\lambda
_{Xj}}\sum_{k}r_{kj}^{2}\biggr]\sum
_{m}\biggl[\sum_{\ell}\langle \varphi
_{\ell},q_{m}\rangle^{2}\biggr]\|R_{YY}v_{m}\|^{2}\\
&&\quad=\sum_{j,k}\frac{r_{kj}^{2}}{\lambda_{Xj}}\sum
_{m}\|R_{YY}v_{m}\|^{2}\qquad \mbox{as }
\|q_{m}\|=1.
\end{eqnarray*}
Note that by \ref{condc2}, the first sum on the right-hand side is bounded. For
the second sum,
\begin{eqnarray*}
\sum_{m}\|R_{YY}v_{m}\|^{2}&=&\sum
_{m}\|R_{YY}^{1/2}q_{m}\|^{2} =\sum_{m}\langle
q_{m},R_{YY}q_{m}\rangle=\sum_{m}\sum_{j}\lambda
_{Yj}\langle q_{m},\varphi_{j}\rangle^{2}\\
&=&\sum_{j}\lambda_{Yj}\sum_{m}\langle q_{m},\varphi
_{j}\rangle^{2} \leq\sum_{j}\lambda_{Yj}\|\varphi
_{j}\|^{2}\leq\sum_{j}\lambda_{Yj}<\infty,
\end{eqnarray*}
which implies (\ref{p1}).
\end{pf*}

\begin{pf*}{Proof of Theorem \protect\ref{th3.4}} Observing
\begin{eqnarray*}
Y_{K}^{\ast}&=&\mathcal{L}_X \beta_{K}^{\ast
}=\sum_{m=1}^{K}\rho
_{m}\mathcal{L}_X(u_{m})R_{YY}v_{m}\\
&=&\sum_{m=1}^{K}\rho
_{m}\langle u_{m},X\rangle R_{YY}v_{m}=\sum_{m=1}^{K}\rho
_{m}U_{m}R_{YY}v_{m}, \\
E\|Y^{\ast}-Y_{K}^{\ast}\|^2&=&E\Biggl\|\sum_{m=K+1}^{\infty
}\rho_{m}U_{m}R_{YY}v_{m}\Biggr\|^{2}=\sum_{m=K+1}^{\infty}\rho
_{m}\|R_{YY}v_{m}\|^{2} \quad \mbox{and}\\
E\|\mathcal{L}_X\beta_{K}^{\ast
}\|^{2}&=&E\Biggl\|\sum_{m=1}^{\infty}\rho
_{m}U_{m}R_{YY}v_{m}\Biggr\|^{2}
\\
&=&\sum_{m,j=1}^{\infty
}\rho_{m}\rho_{j}E[U_{m}U_{j}]\langle R_{YY}v_{m},R_{YY}v_{j}\rangle
=\sum_{m=1}^{\infty}\rho_{m}^{2}\|R_{YY}v_{m}\|^{2}<\infty,
\end{eqnarray*}
we infer that $E\|Y^{\ast}-Y_{K}^{\ast}\|^{2}\rightarrow0$ as
$K\rightarrow\infty.$ From $E[U_{m}]=0,$ for $m\geq1,$ we have
$E[Y_{K}^{\ast}]=0$ and, moreover,
\begin{eqnarray*}
E\|Y-Y_{K}^{\ast}\|^{2}&=&E\|(Y-\mathcal{L}_X \beta_{0}^{\ast
})+\mathcal{L}_X(\beta_{0}^{\ast}-\beta_{K}^{\ast})\|^{2}
\\
&=&E\|Y-\mathcal{L}_X \beta_{0}^{\ast
}\|^{2}+E\|\mathcal{L}_X(\beta_{0}^{\ast}-\beta_{K}^{\ast
})\|^{2}+2E\langle Y-\mathcal{L}_X \beta_{0}^{\ast
},\mathcal{L}_X(\beta_{0}^{\ast}-\beta_{K}^{\ast})\rangle.
\end{eqnarray*}
Since $E\|Y-\mathcal{L}_X \beta_{0}^{\ast
}\|^{2}=\operatorname{trace}(R_{YY})-E\|\mathcal{L}_X\beta_{0}^{\ast
}\|^{2}$ and as $\beta_{0}^{\ast}$ is the solution of the normal
equation (\ref{norm}), we obtain $E\langle Y-\mathcal{L}_X \beta
_{0}^{\ast},\mathcal{L}_X(\beta_{0}^{\ast}-\beta_{K}^{\ast
})\rangle=E\langle \mathcal{L}_X^{\ast}(Y-\mathcal{L}_X \beta
_{0}^{\ast}),\beta_{0}^{\ast}-\beta_{K}^{\ast}\rangle=0.$
Likewise,
\[
E\|\mathcal{L}_X(\beta_{0}^{\ast}-\beta_{K}^{\ast
})\|^{2}=\sum_{m=K+1}^{\infty}\rho
_{m}^{2}\|R_{YY}v_{m}\|^{2},
\]
implying (\ref{conv}).
\end{pf*}

\begin{pf*}{Proof of Theorem \protect\ref{th6.3}} From (\ref{canreg2}),
(\ref{k-pred}) for any $K\geq1$,
\[
R_{Y^{\ast}Y^{\ast
}}-R_{Y_{K}^{\ast}Y_{K}^{\ast
}}=R_{YY}^{1/2}\Biggl[\sum_{m=K+1}^{\infty}\rho
_{m}^{2}q_{m}\otimes q_{m}\Biggr] R_{YY}^{1/2}=R_{YY}^{1/2}R_{K+1}^{\ast
}R_{K+1}R_{YY}^{1/2},
\]
where $R_{K+1}=\operatorname{Proj}_{\operatorname{span}\{q_{m},m\geq K+1\}}R$ and
hence, $R_{Y^{\ast}Y^{\ast}}-R_{Y_{K}^{\ast}Y_{K}^{\ast}}\geq
0.$ Note that
\begin{eqnarray*}
r_{Y_{c,\infty}Y_{c,\infty}}(s,t)&=&E[Y_{c,\infty}(s)Y_{c,\infty
}(t)]=\sum_{m,j=1}^{\infty
}E[V_{m}V_{j}]R_{YY}v_{m}(s)R_{YY}v_{j}(t)
\\
&=&\sum_{m=1}^{\infty
}R_{YY}v_{m}(s)R_{YY}v_{j}(t)=\sum_{m=1}^{\infty
}R_{YY}^{1/2}(q_{m})(s)R_{YY}^{1/2}(q_{m})(t),
\end{eqnarray*}
implying that
\begin{eqnarray*}
R_{Y_{c,\infty}Y_{c,\infty}}-R_{Y^{\ast}Y^{\ast
}}=R_{YY}^{1/2}\Biggl[\sum_{m=1}^{\infty}(1-\rho
_{m}^{2})q_{m}\otimes q_{m}\Biggr]R_{YY}^{1/2}\geq0.
\end{eqnarray*}
Finally, from Lemma \ref{le6.2}(b), we have $Y=Y_{c,\infty}-Y_{c,\infty
}^{\bot},$ therefore $r_{YY}=r_{Y_{c,\infty}Y_{c,\infty
}}+r_{Y_{c,\infty}^{\bot}Y_{c,\infty}^{\bot}}.$ This leads to
$r_{YY}-r_{Y_{c,\infty}Y_{c,\infty}}=r_{Y_{c,\infty}^{\bot
}Y_{c,\infty}^{\bot}}$ and $R_{YY}-R_{Y_{c,\infty}Y_{c,\infty
}}=R_{Y_{c,\infty}^{\bot}Y_{c,\infty}^{\bot}}\geq0.$
\end{pf*}

We need the following auxiliary result to prove Theorem \ref{th6.3}. We call
two $L_{2}$-processes $X$ and~$Y$ \textit{uncorrelated} if and only if
$E[\langle
u,X\rangle\langle v,Y\rangle]=0$ for all $L_{2}$-functions $u$ and $v$.

\begin{lemma}\label{le7.1}$Y_{c,\infty}^{\bot}$ and
$Y^{\ast}$ are uncorrelated.
\end{lemma}

\begin{pf} For any $\tilde{u},$ $\tilde{v}\in
L_{2}(T_{2})$, write
$\tilde{v}=\tilde{v}_{1}+\tilde{v}_{2}$, with
$R_{YY}^{1/2}$\ $\tilde{v}_{1}\in\operatorname{span}\{q_{m}; m\geq
1\}$, which is equivalent to $\tilde{v}_{1}\in
\operatorname{span}\{v_{m}; m\geq1\}$ and $R_{YY}^{1/2}$\ $
\tilde{v}_{2}\in\operatorname{span}\{q_{m}; m\geq1\}^{\bot}$. Then
\[
\langle \tilde{v}_{2},Y^{\ast}\rangle=\sum_{m=1}^{\infty
}\rho_{m}U_{m}\langle
\tilde{v}_{2},R_{YY}v_{m}\rangle=\sum_{m=1}^{\infty}\rho
_{m}U_{m}\langle R_{YY}^{1/2}\tilde{v}_{2},q_{m}\rangle=0.
\]
With $\tilde{v}_{1}=\sum_{m}a_{m}v_{m},$ write $\langle
\tilde{v},Y^{\ast}\rangle=\langle \tilde{v}_{1},Y^{\ast
}\rangle=\sum_{m,j}a_{m}\rho_{j}U_{j}\langle
v_{m},R_{YY}v_{j}\rangle=\sum_{m}a_{m}\rho_{m}U_{m}.$
Furthermore, from Lemma \ref{le6.2}(b), $E[U_{m}\langle
\tilde{u},Y_{c,\infty}^{\bot}\rangle]=0$ for all $m\geq1$. We
conclude that $E[\langle \tilde{u},Y_{c,\infty}^{\bot}\rangle
\langle
\tilde{v},Y^{\ast}\rangle]=0.$
\end{pf}

\begin{pf*}{Proof of Theorem \protect\ref{th6.4}} Calculating the covariance
operators for $(Y,Y^{\ast})$,
\begin{eqnarray*}
r_{Y^{\ast}Y^{\ast}}(s,t)&=&E[Y^{\ast}(s)Y^{\ast}(t)]
=\sum_{m,j}\rho_{m}\rho
_{j}E[U_{m}U_{j}]R_{YY}u_{m}(s)R_{YY}v_{j}(t)
\\
&=&\sum_{m}\rho_{m}^{2}R_{YY}u_{m}(s)R_{YY}v_{m}(t)
=\sum_{m}\rho
_{m}^{2}R_{YY}^{1/2}q_{m}(s)R_{YY}^{1/2}q_{m}(t)
\end{eqnarray*}
so that
\[
R_{Y^{\ast}Y^{\ast}}=\sum_{m}\rho
_{m}^{2}R_{YY}^{1/2}[q_{m}\otimes q_{m}]R_{YY}^{1/2}
=R_{YY}^{1/2}\biggl[\sum_{m}\rho_{m}^{2}q_{m}\otimes
q_{m}\biggr]R_{YY}^{1/2}=R_{YY}^{1/2}R_{0}R_{YY}^{1/2}.
\]
Now, from Lemmas \ref{le6.2} and \ref{le7.1},
\begin{eqnarray*}
r_{YY^{\ast}}(s,t)&=&E[Y(s)Y^{\ast}(t)]=E\bigl[\bigl(Y_{c,\infty
}(s)+Y_{c,\infty}^{\bot}(s)\bigr)Y^{\ast}(t)\bigr]
\\
&=&E[Y_{c,\infty}(s)Y^{\ast
}(t)]=E\biggl[\sum_{m}V_{m}R_{YY}v_{m}(s)\sum_{j}\rho
_{j}U_{j}R_{YY}v_{j}(t)\biggr]
\\
&=&\sum_{m,j}E[V_{m}U_{j}\rho
_{j}R_{YY}v_{m}(s)R_{YY}v_{j}(t)]\\
&=&\sum_{m}\rho
_{m}^{2}R_{YY}v_{m}(s)R_{YY}v_{j}(t)=r_{Y^{\ast}Y^{\ast}}(s,t).
\end{eqnarray*}
Hence, $R_{YY^{\ast}}=R_{Y^{\ast}Y^{\ast}}$. The correlation
operator for $(Y,Y^{\ast})$ is $\tilde{R}=R_{YY}^{-1/2}R_{YY^{\ast
}}R_{Y^{\ast}Y^{\ast}}^{-1/2}=R_{YY}^{-1/2}R_{Y^{\ast}Y^{\ast
}}^{1/2}$ with $\tilde{R}\tilde{R}^{\ast
}=R_{YY}^{-1/2}R_{Y^{\ast}Y^{\ast}}R_{YY}^{-1/2}=R_{0}.$ Hence,
$\tilde{\rho}_{m}=\rho_{m},$ $\tilde{p}_{m}=q_{m}$ and
$\tilde{q}_{m}=\tilde{R}^{\ast}
\tilde{p}_{m}/\tilde{\rho}_{m}=R_{Y^{\ast}Y^{\ast
}}^{1/2}R_{YY}^{-1/2}q_{m}/\rho_{m}=R_{Y^{\ast}Y^{\ast
}}^{1/2}v_{m}/\rho_{m}. $ Moreover,
$\tilde{u}_{m}=R_{YY}^{-1/2}\tilde
{p}_{m}=R_{YY}^{-1/2}q_{m}=v_{m}$ and $\tilde{v}_{m}=R_{Y^{\ast
}Y^{\ast
}}^{-1/2}\tilde{q}_{m}=R_{Y^{\ast}Y^{\ast}}^{-1/2}R_{Y^{\ast
}Y^{\ast}}^{1/2}v_{m}/\rho_{m}=v_{m}/\rho_{m}.$ Note that
$Y_{c,\infty
}^{\ast}=\sum_{m}\tilde{V}_{m}R_{Y^{\ast}Y^{\ast}}
\tilde{v}_{m}$
with
\begin{eqnarray*}
\tilde{V}_{m}&=&\langle \tilde{v}_{m},Y^{\ast}\rangle
=\biggl\langle
v_{m}/\rho_{m},\sum_{j}\rho
_{j}U_{j}R_{YY}v_{j}\biggr\rangle=\sum_{j}U_{j}\langle
v_{m},R_{YY}v_{j}\rangle=U_{m},
\\
R_{Y^{\ast}Y^{\ast
}}\tilde{v}_{m}&=&R_{YY}^{1/2}R_{0}R_{YY}^{1/2}v_{m}/ \rho
_{m}=R_{YY}^{1/2}R_{0}q_{m}/\rho_{m}=\rho
_{m}R_{YY}^{1/2}q_{m}=\rho_{m}R_{YY}v_{m}.
\end{eqnarray*}
Substituting into the equation on the left-hand side of (\ref{t6.4}), one
obtains the equation on the right-hand side of (\ref{t6.4}).
\end{pf*}

\begin{pf*}{Proof of Proposition \protect\ref{pr6.5}} From the definition,
$\mathcal{L}_X^{\ast}$ must satisfy
$\langle \mathcal{L}_X\beta,$ $z\rangle=\langle \beta,\mathcal
{L}_X^{\ast}z\rangle$\break for $%
\beta\in L_{2}(T_{1}\times T_{2})$ and $z\in L_{2}(T_{2}).$ Note
that $\langle \mathcal{L}_X\beta,$
$z\rangle=\int_{T_{2}}(\mathcal{L}_X\beta
)(t)z(t)\,\mathrm{d}t=\break \int_{T_{2}}\int_{T_{1}}X(s)\beta(s,t)z(t)\,\mathrm{d}s\,\mathrm{d}t$ and
$\langle
\beta,\mathcal{L}_X^{\ast}z\rangle=\int\int_{T_{1}\times
T_{2}}\beta(s,t)(\mathcal{L}_X^{\ast}z)(s,t)\,\mathrm{d}s\,\mathrm{d}t.$ For the
differences, we obtain $\int\int\beta
(s,t)[X(s)z(t)-(\mathcal{L}_X^{\ast}z)(s,t)]\,\mathrm{d}s\,\mathrm{d}t=0$ for arbitrary
$\beta\in L_{2}(T_{1}\times T_{2})$ and $z\in L_{2}(T_{2})$. This
implies that $(\mathcal{L}_X^{\ast}z)(s,t)=X(s)z(t).$
\end{pf*}

\begin{pf*}{Proof of Proposition \protect\ref{pr6.6}} By Proposition \ref{pr6.5},
$\Gamma_{XX}=E[\mathcal{L}_X^{\ast}\mathcal{L}_X]$. Since the integral
operator $\Gamma_{XX}$ has the $L_{2}$-integral kernel $r_{XX}$, it
is a Hilbert--Schmidt operator (Conway (\citeyear{Conway85})). Moreover, for $\beta
_{1},\beta_{2}\in L_{2}(T_{1}\times T_{2})$,
\begin{eqnarray*}
\langle \Gamma_{XX}\beta_{1},\beta_{2}\rangle&=&\int\int(\Gamma
_{XX}\beta_{1})(s,t)\beta_{2}(s,t)\,\mathrm{d}s\,\mathrm{d}t =\int\int\int
r_{XX}(s,w)\beta_{1}(w,t)\beta_{2}(s,t)\,\mathrm{d}w\,\mathrm{d}s\,\mathrm{d}t,\\
\langle \beta_{1}, \Gamma_{XX}\beta_{2}\rangle&=&\int\int\beta
_{1}(s,t)(\Gamma_{XX}\beta_{2}(s,t))\,\mathrm{d}s\,\mathrm{d}t =\int\int\int\beta
_{1}(w,t)r_{XX}(s,w)\beta_{2}(s,t)\,\mathrm{d}w\,\mathrm{d}s\,\mathrm{d}t,
\end{eqnarray*}
implying that $\Gamma_{XX}$ is self-adjoint. Furthermore, $\Gamma
_{XX}$ is non-negative definite because, for arbitrary $\beta\in
L_{2}(T_{1}\times T_{2})$,
\begin{eqnarray*}
\langle \Gamma_{XX}\beta,\beta\rangle&=& \int\int\int
E[X(s)X(w)]\beta(w,t)\beta(s,t)\,\mathrm{d}w\,\mathrm{d}s\,\mathrm{d}t\\
&=&E\biggl[\int(\mathcal{L}_X\beta)(t)(\mathcal{L}_X\beta
)(t)\,\mathrm{d}t\biggr]=E\|\mathcal{L}_X\beta\|^{2}\geq0.
\end{eqnarray*}
\upqed\end{pf*}

\section*{Acknowledgements}
We wish to thank two referees for careful reading and are especially
indebted to one reviewer and the Associate Editor for comments which
led to substantial changes and various corrections. This research
was supported in part by NSF Grants DMS-03-54448, DMS-04-06430,
DMS-05-05537 and DMS-08-06199.

\printhistory

\end{document}